\documentclass[letterpaper, 10 pt, conference]{ieeeconf}  

\IEEEoverridecommandlockouts                              

\usepackage{amsfonts}
\usepackage{amsmath, amssymb}
\usepackage[utf8]{inputenc}
 \usepackage{amsfonts}
\usepackage{mathtools}
\usepackage{enumerate}
\usepackage{tikz}
\usepackage{subfigure}
\usepackage{graphics,color}
\usepackage[noend]{algpseudocode}
\usepackage{algorithm}
\usepackage{algorithmicx}

\newcommand{\R}{\mathbb{R}}

\def\balign{\begin{align*}}
\def\ealign{\end{align*}}
\def\bbmat{\begin{bmatrix}}
\def\ebmat{\end{bmatrix}}

\newcommand*\Let[2]{\State #1 $\gets$ #2}

\newtheorem{theorem}{Theorem}
\newtheorem{example}{Example}
\newtheorem{lemma}{Lemma}
\newtheorem{corollary}{Corollary}

\newtheorem{proposition}{Proposition}
\newtheorem{assumption}{Assumption}
\newtheorem{definition}{Definition}

\newtheorem{remark}{Remark}
\def\BT{\begin{theorem}}
\def\ET{\end{theorem}}
\def\BL{\begin{lemma}}
\def\EL{\end{lemma}}
\def\BP{\begin{proposition}}
\def\EP{\end{proposition}}
\def\BC{\begin{corollary}}
\def\EC{\end{corollary}}
\def\BD{\begin{definition}}
\def\ED{\end{definition}}
\def\BA{\begin{assumption}}
\def\EA{\end{assumption}}
\def\BR{\begin{remark}}
\def\ER{\end{remark}}
\def\BE{\begin{example}}
\def\EE{\end{example}}

\title{\LARGE \bf
A Relaxed Optimization Approach for Cardinality-Constrained Portfolio Optimization 
}

\author{Jize Zhang $^{1}$, Tim Leung $^{2}$, Aleksandr Aravkin $^{3}$
\thanks{*This work was supported by the Washington Research Foundation Data Science Professorship.}
\thanks{$^{1,2,3}$ Department of Applied Mathematics, University of Washington, Seattle.
        {\tt\small jizez@uw.edu, timleung@uw.edu, saravkin@uw.edu}}%
}

\begin{document}

\maketitle
\thispagestyle{empty}
\pagestyle{empty}

\begin{abstract}
A cardinality-constrained portfolio caps the number of stocks to be traded across and within groups or sectors.  These limitations arise from real-world scenarios faced by fund managers, who are constrained by transaction costs and client preferences as they seek to  maximize return and limit risk. 

We develop a new approach to solve cardinality-constrained portfolio optimization problems, 
extending both Markowitz and conditional value at risk (CVaR) optimization models with cardinality constraints. We derive a continuous relaxation method for the NP-hard objective, 
which allows for very efficient algorithms with standard convergence guarantees for nonconvex problems.  For smaller cases, where brute force search is feasible to compute the globally optimal cardinality-constrained portfolio,  the new approach finds the best portfolio for the cardinality-constrained Markowitz model  and a very good local minimum for the cardinality-constrained CVaR model. For higher dimensions, where brute-force search 
is prohibitively expensive, we find feasible portfolios that are nearly as efficient 
as their non-cardinality constrained counterparts.  

\end{abstract}

\section{INTRODUCTION}

 {Since the introduction of the original Markowitz \cite{markowitz} mean-variance portfolio selection model, 
there have been many extensions and refinements by both practitioners and academics. 
Among them, one type of constraints arises from industry practice. 
For many portfolio managers, it is desirable or even imperative to limit the number of assets in a portfolio 
and/or impose limits on the proportion of the portfolio devoted to any particular asset or asset class. 
These constraints can be driven by both the portfolio mandates set by clients (investors), 
or pragmatic reasons such as transaction costs, minimum lot sizes, and execution efficiency. }

Prior art for cardinality constrained portfolio selection comprises two classes of methods: heuristic algorithms such as genetic algorithms and particle swarm algorithms~\cite{chang2000heuristics,soleimani2009markowitz,deng2012markowitz}, 
and mixed integer programming~\cite{shaw2008lagrangian,bertsimas2009algorithm,Cui2013}.
Both lines of research have primarily focused on the Markowitz (mean-variance) criterion. 
One exception is~\cite{murray2012local}, which projects asset returns onto a reduced space, then 
uses clustering to identify similar subgroups;~\cite{murray2012local} also assumes a quadratic loss function.


In this paper, we propose a different approach by formulating the problem as a continuous optimization problem over a highly nonconvex set induced by the intersection 
of cardinality and simplex constraints.  We then develop a relaxation method using auxiliary variables, and create  an efficient projection map onto the nonconvex set.  These innovations allow recently developed techniques for  structured nonsmooth nonconvex optimization~\cite{bolte2014proximal} to bear on the problem. 
The  proposed approach is not limited to Markowitz portfolio selection, but it  can also be applied to a variety  of portfolio criteria, both smooth and nonsmooth. 
As a key example, we look at the nonsmooth CVaR portfolio selection problem with cardinality constraints.

The paper is organized as follows. In Section \ref{sec:port_models} we review Markowitz and CVaR portfolio optimization formulations.
In Section \ref{sec:constr} we present the cardinality constraints, along with the relaxation method that leads to an efficient solution to the constrained problem. Section \ref{sec:proj} develops an efficient projection 
onto the nonconvex set induced by cardinality and simplex constraints, and presents a proof of correctness of this projection. 
Section \ref{sec:algorithm} then discusses general algorithm and its specific instantiations  for Markowitz and CVaR portfolio selection models, 
and characterizes the stationarity points of the relaxed problem.  
 Section \ref{sec:numeric} shows numerical results using a dataset with 65 stocks, including  a comparison with brute force search for very small selections, 
  and a larger-scale experiment that compares efficient frontiers of cardinality-constrained and standard portfolios. 

\section{Portfolio Selection Models}
\label{sec:port_models}
The general portfolio selection problem can be formulated as follows. 
We are given a total of $n$ candidate assets and a certain selection criterion, usually formulated through an objective $f(w)$. 
Portfolio selection models also impose bounds on $\omega$,  often using the simplex:
\[
\Delta_1(w) = \left\{w: 0 \leq w_i \leq 1, \; 1^Tw = 1\right\}.
\]
The basic portfolio optimization problem is given by 
\begin{equation}
\label{eq:port_simple}
\min_{w \in \Delta_1} f(w) 
\end{equation}
with two widely used objectives $f(w)$ given below. 

\begin{example}[Mean-Variance]
In mean-variance portfolio selection, also known as Markowitz selection, 
\begin{equation}
\label{eq:objMarkowitz} 
f(w) = w^T\Sigma w - \gamma \mu^T w
\end{equation}
where $\Sigma \in \R^{n \times n}$ is the covariance matrix  
and $\mu \in \R^n$ the expected return vector. 
Quantities $\Sigma$ and $\mu$ are computed from historical returns:
\[
\mu = \frac{1}{N}\sum_{j=1}^N r_j, \quad \Sigma = \frac{1}{N}\sum_{i=1}^N (r_j - \mu) (r_j - \mu)^T,
\]
where $r_j$ are historical return vectors, and $N$ is the total number of samples. 
The parameter $\gamma$ controls the weight between variance (as a measure of risk) and return. 
\end{example} 

\begin{example}[Conditional Value-at-Risk]
The conditional value-at-risk (CVaR) model~\cite{rockafellar2002conditional} minimizes the CVaR superquantile 
over portfolio selections: 
\begin{equation}
\label{eq:objCVaR}
f_\beta(w) = \min_\alpha \alpha + \frac{1}{N(1-\beta)}\sum_{j=1}^N [-w^Tr_j - \alpha]_+, 
\end{equation}
where the quantile $\beta$ is related to the CVaR value $\alpha$  by  
\[
P(\mbox{loss} > \alpha) \leq 1-\beta.
\]
\end{example}

\section{Combinatorial Constraints}
\label{sec:constr}

Combinatorial constraints restrict the number of stocks to purchase, within specified subgroups and/or across the entire portfolio. 
We consider the constraint set $\Omega$ given by 
\begin{equation}
\label{eq:omega}
\Omega= \left\{p_i \leq 1^Tw_i \leq q_i, \; \|w_i\|_0 \leq k_i \in \mathbb{N}_+, \; i = 1, .., m\right\}, 
\end{equation}
where the portfolio weights $w$ are divided into subgroups $w = [w_1,...,w_i,...,w_m]$, and 
\begin{equation}
\|w_i\|_0 = \mbox{card}(j: w_i^j \neq 0).\label{l0norm}
\end{equation}
 The $\ell_0$ norm in \eqref{l0norm} counts the number of nonzero entries of its argument. 
The set $\Omega$ captures a wide variety of realistic trading constraints. 
We can group assets by sectors, industries or other criteria. The constraints restrict both the fraction of the portfolio, as well as 
certain number of assets from each group.
For example, we can require that the trader buy no more than 5 stocks from healthcare, 
between 3 and 8 stocks from tech, and that stocks in consumer goods should comprise between 10\% to 25\% of the 
total portfolio.

We now consider the constrained optimization problem 
\begin{equation}
\label{eq:portOpt}
\min_{w \in \Omega \cap \Delta_1} f(w) 
\end{equation}
with $\Omega$ as in~\eqref{eq:omega},  simplex constraint $\Delta_1$, 
and portfolio criterion $f(w)$, such as Markowitz~\eqref{eq:objMarkowitz} or CVaR~\eqref{eq:objCVaR}.
The inclusion of the $l_0$-norm in $\Omega$ means that $w$ lies in the intersection of a highly nonconvex set and 
a compact convex set given by box-constraints and the 1-simplex. 

To develop an efficient method for~\eqref{eq:portOpt}, we first relax the problem by introducing an auxiliary variable $v \in \R^m$ 
to lie in $\Omega$, and forced to be close to $w \in \Delta_1$ using a quadratic penalty term. 
This type of relaxation was recently shown to be effective for nonsmooth nonconvex optimization~\cite{zheng2018fast}.
The relaxed problem is given by 
\begin{equation}
\label{eq:port_relaxed}
\begin{aligned}
 \min_{w, v} & \quad f(w) + \frac{\nu}{2}\|w-v\|^2\\ 
\mbox{s.t.} &\quad v \in  \Omega, \quad w \in \Delta_1 
\end{aligned}
\end{equation}
As $\nu$ increases, we have $\|w - v\| \leq \frac{C}{\nu}$ for some constant $C$, and problem~\eqref{eq:port_relaxed}
approximates~\eqref{eq:portOpt}.

\subsection{$l_1$ norm vs $l_0$ norm}
Since $l_0$ norm constraints are highly nonconvex in general, it is common to replace them with convex $l_1$ norm constraints. 
However, for portfolio optimization 1-norm constraints are not a useful relaxation, since $w \in \Delta_1$ forces $\|w\|_1 = 1$. 
Thus requiring $\|w\|_1 \leq \tau$ makes the problem infeasible for $\tau < 1$ and is meaningless for $\tau > 1$.
For cardinality-constrained optimization, we need the $l_0$ norm for its simple and direct interpretation: 
$\|w\|_0 \leq k$ means exactly that we allow no more than $k$ assets. 

Even though the constraint is very nonconvex, there is a simple form for the projection onto the set $\mathbb{B}_{0}^k := \|w\|_0 \leq k$. 
Let $\omega_k(z)$ be the set of indices corresponding to the largest $k$ entries of $z$ (by absolute value). Then  
\[
\mbox{Proj}_{\mathbb{B}_0^k}(z)_i = \begin{cases}
z_i & \quad i \in \omega_k(z) \\
0 & \quad i \not\in \omega_k(z). 
\end{cases} 
\]
The projection onto the entire set $\Omega$ is more complicated, and is developed in the next section. 



\section{Projection Map}
\label{sec:proj}
In this section we develop the projection onto $\Omega$ as defined in~\eqref{eq:omega}. 
First we introduce a useful lemma for this projection.
The proofs of all technical lemmas are given in the Appendix. 
\begin{lemma}
\label{lm:proj}
Suppose $y \in \mathbb{R}^{l}$. Let $K$ be any size-$k$ subset of $I = \{1,...,l\}$ and $\mathcal{K}$ the union of all such $K$s, so that $I-K$ denotes the complement of $K$ in $I$. Let $\mathcal{C} \in \R_+$ be a closed convex subset of $\R_+ = \{ x: x \geq 0\}$. Without loss of generality, we reorder $y$ such that $y_1 \geq y_2 \geq ...\geq y_l$. We claim that the optimal $K$ for the problem 
\[ \min_{z_K \in \mathcal{C},z_{I-K} = 0,K \in \mathcal{K}} \frac{1}{2}\|y-z\|^2 \] 
is $K_{opt} = \{1,2,..,k\}$, i.e. the indices corresponding to the $k$ largest components in $y$.
\end{lemma}

Applying Lemma \ref{lm:proj} to the projection of each subgroup $i$ in~\eqref{eq:omega}, 
for each $w_i$ we pick its $k_i$-largest components and project them onto the set 
\[
\left\{ z \geq 0, \quad p_i \leq 1^Tz \leq q_i\right\}.
\]
The problem thus reduces to solving $m$ independent projection subproblems  
\begin{equation}
\label{eq: proj_sub}
\min_{z \geq 0, \quad p_i \leq 1^Tz \leq q_i} \|u_i-z\|^2, \quad v \in \R^{k_i} 
\end{equation}
where $u_i\in\mathbb{R}^k$ contains any selection of $k_i$-largest components in $w_i$. 
This is a quadratic problem with linear constraints and can be solved exactly. In particular we relate it to projection onto simplex via the following lemma. 

\begin{lemma}
\label{lem:box}
Consider the problem
\[\min_{z \geq 0, p\leq 1^Tz \leq q} \|u-z\|^2, z \in \R^{k}, \quad 0 \leq p \leq q \leq 1.\] 
Let $m$ be the number of positive entries in $u$ and $z^*$ the minimizer. Then we have 
\begin{eqnarray}
1^Tu_{1:m} < p & \Rightarrow & z^* = \text{argmin}_{z \in \Delta_{p}}\|u-z\|^2 \\
p \leq 1^Tu_{1:m} \leq q & \Rightarrow & z^*_{1:m} = u_{1:m}, z^*_{m+1:k} = 0 \\
1^Tu_{1:m} > q & \Rightarrow & z^* = \text{argmin}_{z \in \Delta_q}\|u-z\|^2.
\end{eqnarray}
\end{lemma}

\section{Algorithms}
\label{sec:algorithm}
To solve Problem (\ref{eq:port_relaxed}) we use proximal alternating linearized minimization (PALM) \cite{bolte2014proximal}, 
with alternating updates on $w$ and $v$. 
The update on $v$ is always a projection, whereas the update on $w$ is problem dependent. 
%

\subsection{Cardinality-Constrained Markowitz}
The relaxed problem for mean-variance (Markowitz) portfolio optimization is 
\[ 
\begin{aligned}
\min_{v,w \in \Delta_1} g(v,w):= &w^T(\Sigma + \lambda I) w -\gamma \mu^Tw + \sum_{i=1}^m \frac{\nu}{2}\|w_i-v_i\|^2 \\
\mbox{s.t.} & \;  v_i \geq 0, p_i \leq 1^Tv_i \leq q_i, \|v_i\|_0 \leq k_i.
\end{aligned}
\]
where $\lambda I$ is a ridge regularization term in case $\Sigma$ is singular. Since the problem involves only one nonsmooth term in $w$ (the simplex constraint), we can use proximal gradient step directly to update $w$, as shown in Algorithm~\ref{alg:mvar}.
\begin{algorithm}
  \caption[Caption]{\label{alg:mvar} PALM for Cardinality-Constrained Markowitz} 
  \begin{algorithmic}[1]
    \Require{$w,v \in \R^n, \gamma,\lambda,\nu$} \\
    $\mathcal{V}_i = \{ z \geq 0: p_i \leq 1^Tz \leq q_i ,\|z\|_0 \leq k_i\}$
   \While{not converged} 
      	\Let{$w$}{Proj$_{\Delta_1}(w - \delta (\nabla_w f(w) + \nu(w-v)) $} 
	\For{$i = 1,...,m$}
	\Let{$v_i$}{Proj$_{\mathcal{V}_i} (w_i)$}
	    \EndFor
	 \EndWhile
      \Return{$w,v$}
  \end{algorithmic}
\end{algorithm}

\subsection{Cardinality-Constrained CVaR}
The relaxed objective for CVar is 
\begin{align*}
\min_{\alpha,w \in \Delta_1,v}& \quad g(w,\alpha,v )  \\
\mbox{s.t.} & \quad v_i \geq 0, p_i \leq 1^Tv_i \leq q_i, \|v_i\|_0 \leq k_i,
\end{align*}
where
 \[ 
 g(w,\alpha,v ) = \alpha + \frac{1}{N(1-\beta)}\sum_{j=1}^N [-w^Tr_j - \alpha]_+ + \sum_{i=1}^m \frac{\nu}{2}\|w_i-v_i\|^2 .\]
Unlike the mean variance problem, here the problem has two nonsmooth terms in $w$, the simplex constraint and the 
hinge loss $[-w^Tr_j - \alpha]_+$, which we also relax. 
We introduce a second auxiliary variable $u \in \R^N$ and let $-Rw -\alpha = u$, 
$R \in \R^{N \times n}$ is the matrix of all samples $r_j$ stacked together. The objective now becomes
\begin{align*}
 \min_{\alpha,w \in \Delta_1,v,u} & \quad \alpha + \frac{1}{N(1-\beta)}\sum_{j=1}^N [u_j]_+   \\
 & \quad + \frac{\gamma}{2}\|Rw+\alpha{\bf1}+ u\|^2 + \sum_{i=1}^m \frac{\nu}{2}\|w_i-v_i\|^2 \\
\mbox{s.t.} & \quad v_i \geq 0, p_i \leq 1^Tv_i \leq q_i, \|v_i\|_0 \leq k_i.
 \end{align*}
We can partially minimize in $\alpha$, since the objective with respect to $\alpha$ is an unconstrained quadratic: 
\[
 \alpha^*(u,w) = -\frac{1 + \gamma {\bf 1}^T(Rw + u)}{\gamma N}.
 \]
Plugging in $\alpha^*$, the problem simplifies to 
\begin{equation}
\label{eq:cvar2}
\begin{aligned}
\min_{w \in \Delta_1, u, v \geq 0}  &\quad  \widetilde g(w,u,v) \\
\mbox{s.t.} & \quad v_i \geq 0, p_i \leq 1^Tv_i \leq q_i, \|v_i\|_0 \leq k_i
\end{aligned}
\end{equation}
where
\begin{equation}
\label{eq:cvarloss}
\begin{aligned}
\widetilde g(w,u,v) =  &\alpha^*(u,w) +  \frac{1}{N(1-\beta)}\sum_{j=1}^N [u_j]_+ + \frac{\nu}{2}\|w-v\|^2 \\
& + \frac{\gamma}{2}\left\|\left(I-\frac{{\bf 1}{\bf 1}^T}{N}\right)(Rw+u) - \frac{1}{\gamma N} {\bf 1}\right\|^2.
\end{aligned}
\end{equation}
Again we use proximal gradient step to update $w$ and $u$. 

\subsection{Acceleration Strategies}
Empirically replacing proximal gradient step for $w$ and $u$ with momentum-based updates  can help speeding up the convergence 
(see Algorithm \ref{alg:fista}). The convergence theory is available only for convex models~\cite{beck2009fast}, 
but FISTA updates are known to help in a lot of nonconvex models.  
The algorithm for~\eqref{eq:cvar2} 
accelerated with FISTA updates is detailed in Algorithm (\ref{alg:cvar2}).
\begin{algorithm}
  \caption[Caption]{\label{alg:fista} FISTA-Update} 
  \begin{algorithmic}[4]
    \Require{ $t_k, y_k, x_k$}
    \Let{$x_{k+1}$}{ ProxGrad($y_k$)}
    \Let{$t_{k+1}$}{$\frac{1+\sqrt{1+4t_k^2}}{2}$}
    \Let{$y_{k+1}$}{$x_k + \frac{t_k-1}{t_{k+1}}(x_k-x_{k+1})$}\\\
\Return{$(x_{k+1}, y_{k+1}, t_{k+1})$}
  \end{algorithmic}
\end{algorithm}

In each iteration we perform one F-update for $w$ and then one for $u$, 
followed by a projection step on $v$. The $y$ and $t$ values corresponding to both $w$ and $v$ 
are saved and used in the next iteration.

\begin{figure}
\centering
\includegraphics[scale = .5]{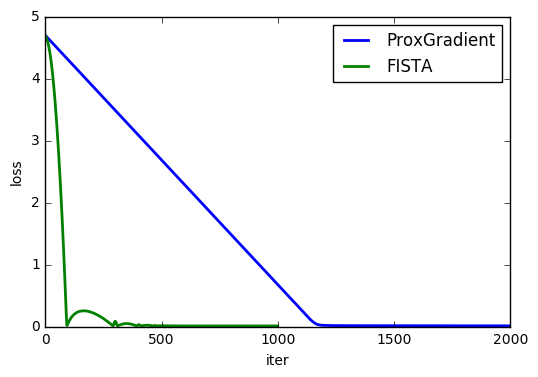}
\caption{FISTA vs prox-gradient for CVaR loss $\widetilde g$~\eqref{eq:cvarloss}, with $\beta = 0.9$.}
\label{fig:cvar_fista}
\end{figure}

Figure \ref{fig:cvar_fista} shows a comparison of FISTA updates and prox-gradient updates on $w$ and $u$. The value of $\beta$ is set to be $0.9$ and the maximum number of stocks selected from each sector is restricted to be 5. Varying the $\beta$ value and constraints yields similar plots.

\begin{algorithm}
  \caption[Caption]{\label{alg:cvar2} Accelerated PALM for CVaR~\eqref{eq:cvar2}} 
  \begin{algorithmic}[3]
    \Require{$w,v \in \R^n, \gamma$} \\
    $\mathcal{V}_i = \{ z \geq 0: p_i \leq 1^Tz \leq q_i ,\|z\|_0 \leq k_i\}$
       \While{not converged} 
     	\Let{$w$}{ \mbox{FISTA-Update}($w;\gamma,\nu$)}
	\Let{$ u $}{ \mbox{FISTA-Update}($u;\gamma,\nu$)}
	\For{$i = 1,...,m$}
	\Let{$v_i$}{Proj$_{\mathcal{V}_i} (w_i)$}
	\EndFor
	\EndWhile
  	\Return{$(w, v, u)$}
  \end{algorithmic}
\end{algorithm}

When there is no sector grouping, the problem reduces to
\[ \min_{w \in \Delta_1}f(w), \|w\|_0 \leq k, \]
we do not need relaxation and can solve the problem using proximal gradient descent. 

\begin{remark}
In some cases, using a continuation on $\nu$ can help avoid local optima. The continuation 
strategy adds an outer loop, increasing $\nu$ as the algorithms proceed. 
\end{remark}

\section{Stationarity}
\label{sec:station}
The PALM algorithm converges to stationary points in nonconvex setting~\cite{bolte2014proximal}.
 A natural question to ask is what stationary points mean in the context of $l_0$-norm constraints. 
In this section, we consider a simple formulation 
\begin{equation}
\label{eq:obj_sta}
\min_{w \in \Delta_1^n} f(w),\quad \|w\|_0  \leq k
\end{equation}
where $f$ is a smooth loss function and $\Delta_1^n$ denotes the 1-simplex in $\R^n$. 
This problem can be solved with proximal gradient descent as discussed in previous section.

We address two questions related to stationarity:
\begin{itemize}
\item[Q$_1$] If we pick $K$ of size $k$, $K \subset \{1,...,n\}$ and solve 
\begin{equation}
\label{eq:obj_sta2}
w^* \in \arg\min_{w} f(w), \quad w_K \in \Delta_1^k, w_j = 0 , j \not\in K
\end{equation}
is $w^*$ automatically a stationary point for~\eqref{eq:obj_sta}?
\item[Q$_2$] Is a stationary point $w^*$ of~\eqref{eq:obj_sta} a fixed point of the proximal gradient method,
as in the convex case? Or is it possible that the algorithm can move to a new  
stationary point with a lower objective value? 
\end{itemize}

Q$_1$ asks whether $w$ being a stationary point of (\ref{eq:obj_sta2}) implies that $w$ is also a stationary point for (\ref{eq:obj_sta}). To answer that question we look at stationarity conditions for both problems. 

Let 
\[
\begin{aligned}
B &\equiv \{ w: w_K \in \Delta_1^k, w_j = 0 , j \not\in K\} \\
 C & \equiv \{ w \in \Delta_1^n: \|w\|_0  \leq k\}. 
\end{aligned}
\]
 Note that $B$ is convex while $C$ is nonconvex. From elementary convex analysis, if $w$ is a stationary point of (\ref{eq:obj_sta2}), then
\[ 0 \in \nabla f(w) + \partial \delta_B(w)\]
where $ \partial \delta_B(w)$ denotes the subdifferential of the indicator function $\delta_B$ at $w$. Since $B$ is convex, $ \partial \delta_B(w)$ is the normal cone to $B$ at $w$, 
defined by
\[ N_B(\bar{w}) = \{ v: \langle v, w - \bar{w}\rangle \leq 0 ~\forall w \in B\}.\]
Similarly if $w$ is a stationary point of (\ref{eq:obj_sta}), then
\[ 0 \in \nabla f(w) + \partial \delta_C(w)\]
where $\partial \delta_C(w)$ is the limiting subdifferential of indicator function $\delta_C$ at $w$. For nonconvex functions, limiting subdifferential is defined in terms of Frechet subdifferentials; see the Appendix (and \cite{rockafellar2009variational} for more details).

We state in Lemma \ref{lm:station} a condition under which a stationary point of (\ref{eq:obj_sta2}) is also a stationary point of (\ref{eq:obj_sta}).
\begin{lemma}
\label{lm:station}
If $w \in \text{int} \Delta_1^k$, then $\partial \delta_C(w) = \partial \delta_B(w)$. This implies that
 $0 \in \nabla f(w) + \partial \delta_B(w) \Leftrightarrow 0 \in \nabla f(w) + \partial \delta_C(w).$
\end{lemma}

\begin{figure}
\centering
\includegraphics[scale = .3]{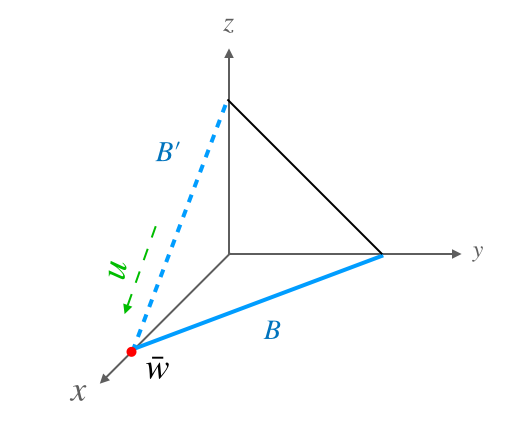}
\label{fig:statn3d}
\caption{
Stationary points for~\eqref{eq:obj_sta} and~\eqref{eq:obj_sta2} only agree on the interior of $\Delta_1^k$; otherwise the former set is smaller. 
To illustrate, take $w \in \R^3$ and require that $\|w\|_0 \leq 2$. Let the line segment in $xy$ plane be $B$ and its left endpoint $\bar{w}$ (red dot). 
Then as $u \to \bar{w}$, it can approach along $B'$ in addition to $B$. In that case there exists $u \in B'$ such that $\langle v, u-\bar{w}\rangle \geq 0$ for some $v \in N_B(\bar{w})$. Hence $v \in N_B(\bar{w})$ but $v \not \in \hat{\partial} \delta_C(\bar{w}).$
}
\end{figure}

Lemma~\ref{lm:station} says that if a stationary point $w$ of (\ref{eq:obj_sta2}) lies in the interior of $\Delta_1^k$, 
then $w$ is also a stationary point of (\ref{eq:obj_sta}). When a $w$ stationary for~\eqref{eq:obj_sta2} is on the boundary, 
it may easily fail to be stationary for~\eqref{eq:obj_sta};
Figure~\ref{fig:statn3d} provides a detailed illustration in $\mathbb{R}^3$. 

Now suppose we solved (\ref{eq:obj_sta2}) and found a stationary point $w$ in the interior of $\Delta_1^k$. 
If we then use the obtained $w$ as an initial guess and solve (\ref{eq:obj_sta}) with proximal gradient descent, would we be trapped at $w$ given it is stationary? 
In other words, is a stationary point necessarily a fixed point for proximal gradient descent? 
The answer is true for convex functions, but not for nonconvex functions. In the convex case, 
we have the following Lemma.

\begin{lemma}
\label{lem:proxcvx}
Consider the problem of 
\[ \min_x f(x) + g(x)\]
where $f(x)$ is smooth and $g(x)$ is convex but nonsmooth. Then for any $\alpha > 0$,  
\[ 0 \in \nabla f(x) + \partial g(x) \Leftrightarrow x = \text{prox}_{\alpha g}(x - \alpha \nabla f(x) ).\]
\end{lemma}

 When $g$ is not convex, e.g. when $g = \delta_C$, then $(I + \alpha \partial g)^{-1}$ can be multi-valued, which means that $ x \in \text{prox}_{\alpha g}(x - \alpha \nabla f(x))$. In that case $x$ is not necessarily a fixed point and we may move to some other points with smaller objective values. Empirically we see that `bad' stationary points are seldom fixed points, 
 and the algorithm will move past these points to stationary points with lower objective values. 
 
 To demonstrate the above claim, we conducted an experiment with the following procedure: 
\begin{enumerate}
\item Randomly pick a subset $K$ of size $k$ from $\{1,...,n\}$
\item Solve the problem $\min_{w} f(w), w_K \in \Delta_1^k, w_j = 0 , j \not\in K$ and denote the minimizer as $w_{init}$
\item Use $w_{init}$ as an initial guess for \eqref{eq:obj_sta} and run our algorithm
\item Check if we find another $w$ with lower objective value
\end{enumerate}
We repeat this procedure with varying $k$ and $n$. For each $k,n$ pair we run 100 trials and record the percentage of times we find a lower objective value;
the results are presented in Table~\ref{table:sim}. 
As $n$ increase we are more likely to stay at $w_{init}$, while as $k$ increases we are more likely to move to a point with a lower objective value. 
When $k=1$ we always stay at $w_{init}$, but in this case we can simply do a linear scan over all assets to find the best one. 
For the mean-variance model, the percentage is high in all cases except $k=1$. 
The study suggests that for moderate $k$ we can expect the algorithm to easily move past a spurious stationary point.  
\begin{table}
\centering
\begin{tabular}{|c|c|c|c|c|}
\hline
Markowitz (n/k) & 1  & 2& 3& 4\\
\hline
15 & 0 & 0.75 & 0.9 & 1  \\
30 & 0 & 0.72 & 0.92 & 0.98\\
45 & 0 & 0.71 & 0.85 & 0.95 \\
\hline \hline
CVaR (n/k) & 1 & 2 &3 &4 \\
\hline 
15 & 0 & 0.31 &0.67  & 0.8 \\
30 & 0& 0.3 & 0.56 & 0.75 \\
45 & 0& 0.2 & 0.42 & 0.64\\
\hline
\end{tabular}
\caption{\label{table:sim} Fraction of times  algorithm decreases objective when starting at a stationary point. }
\end{table}

\section{Numerical Results}
\label{sec:numeric}
In this section we present results from two sets of experiments. The first set of experiments compares our approach against the brute force searches on small examples where a brute force search is feasible. Since gradient based methods only guarantees convergence to local optima for nonconvex problems, this experiment gives a sense of how good 
those local optima are, at least for small cases. Brute-force search is impossible for even moderate-sized problems. 

{The second set of experiments compares Pareto efficient frontiers of models with and without cardinality constraints. 
For Markowitz models, efficient frontiers show the trade-off between portfolio variance against average return; 
for CVaR we plot $\beta$-VaR against $\beta$-CVaR values.}

For all experiments, the underlying dataset comprises closing stock prices for 65 stocks from June 21st 2017 to June 21st 2018, taken from Yahoo finance. 
The stocks belong to 7 sectors, basic metals, consumer goods, finance, health care, industrial goods, service and technology.

\subsection{Comparison with brute force search}
\begin{figure*}
\centering
\begin{tabular}{ccc}
\includegraphics[scale = .35]{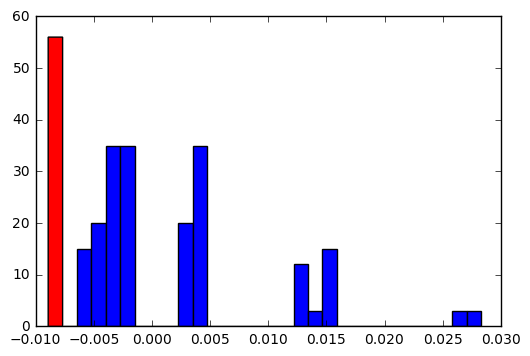} & \includegraphics[scale = .35]{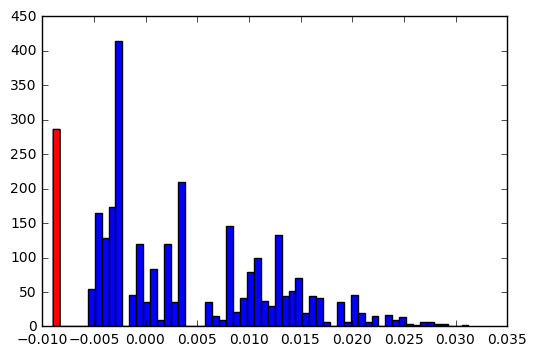}  & \includegraphics[scale = .35]{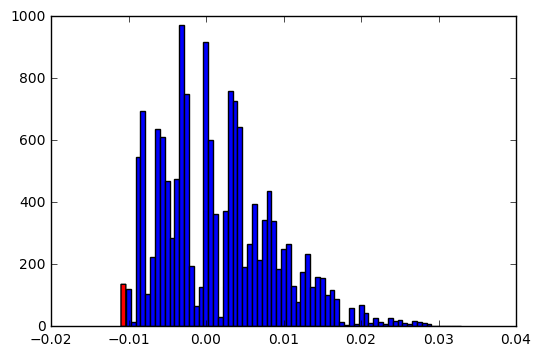} \\
\includegraphics[scale = .35]{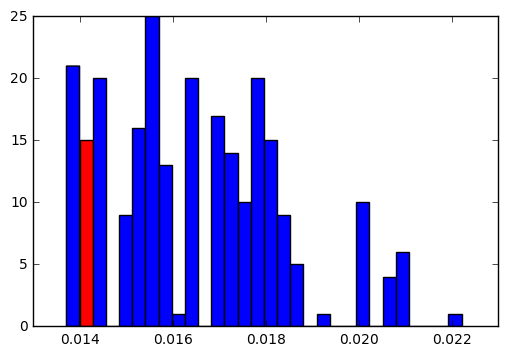} & \includegraphics[scale = .35]{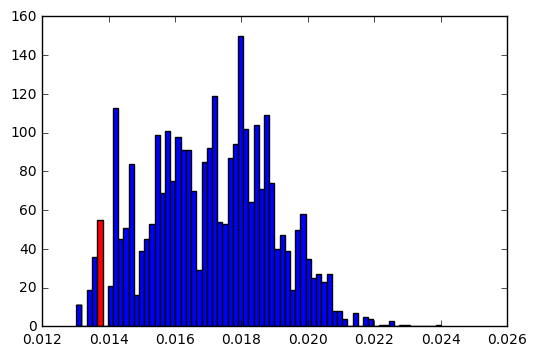} & \includegraphics[scale = .35]{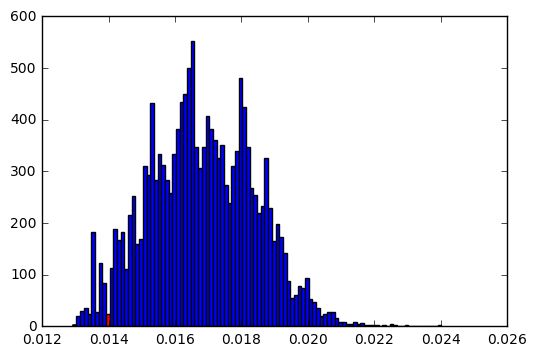}
\end{tabular}
\caption{Histograms of losses. Each plot shows the histogram of loss obtained via exhaustive search and red bar indicates where result from our method lies. Top row: mean-variance; Bottom row: CVaR. Left: 10 choose 5.  Middle: 15 choose 5; Right: 20 choose 5. For mean-variance we found global minima in all three cases; for CVaR we found local minima that are close to global minima.}
\label{cmp_bf}
\end{figure*}

Since the problem is highly nonconvex, our method is not guaranteed to find a global minimum. 
For small simple examples, where exhaustive search is feasible, 
we can compare the quality of our solution (using function value) to that of the global minimum. 
Specifically, we do not impose sector groupings on $w$ (since there is no straightforward way to implement a brute force search in this case) 
and limit the total number of candidate assets. 
Figure \ref{cmp_bf} summarizes the results. The top row shows the histograms of losses for mean-variance portfolio (with $\gamma = 0.1$) 
when we look for at most 5 from a total of 10, 15 or 20 assets (left, middle and right panels). 
The red bar shows the objective value of our solution. The bottom row shows the corresponding histograms for CVaR with $\beta = 0.9$. 
In the experiment, {\bf our solution was a true global minimum in all three simulations for mean-variance portfolios}, and the CVaR solutions were nearly 
optimal in function value. 

As the total number of assets grows, exhaustive search quickly becomes prohibitively expensive. 
For instance, choosing 10 assets out of 30 requires solving more than 30 million optimization problems over the subsets. 
To test our method for these larger scenarios, we first run our algorithm to obtain a solution; 
then we randomly choose asset combinations, solve a minimization problem over that combination, 
and repeat until we find a solution with same or smaller function value. 
We compare the average elapsed time using our method versus brute force search. 
Table \ref{bf_time} shows the results. The task is to choose 10 out of 30 assets and the average is taken over 20 runs. 
Overall, the randomized brute force search takes a very long time to match the quality of the solution (measured by function value) 
that is quickly discovered by continuous optimization. 
 
\begin{table}
\centering
\begin{tabular}{|r|c|c|c|}
\hline
Markowitz & avg time & min time & max time \\
\hline
Proposed method & 0.022 & 0.020 & 0.030 \\
Randomized brute force & 1.77 & 0.018 & 17.39 \\
\hline\hline
CVaR & avg time & min time & max time \\
\hline
Proposed method & 2.03 & 1.80 & 2.64 \\
Randomized brute force & 71.08 & 1.83 & 222.69 \\
\hline
\end{tabular}
\caption{Proposed method vs. randomized brute force search: time (seconds) for choosing 10 assets from 30 (across 20 trials).}
\label{bf_time}
\end{table}

\subsection{Efficient Frontiers}
We plot the Pareto efficient frontier for  classic portfolio optimization strategies against those with different stock cardinalities specified by the cardinality constraints. 
For Markowitz models, an efficient frontier is traced out by varying $\gamma$, the weight on average return, while for CVaR it is traced out by varying the probability level $\beta$.

\begin{figure}
\centering
\includegraphics[scale = .35]{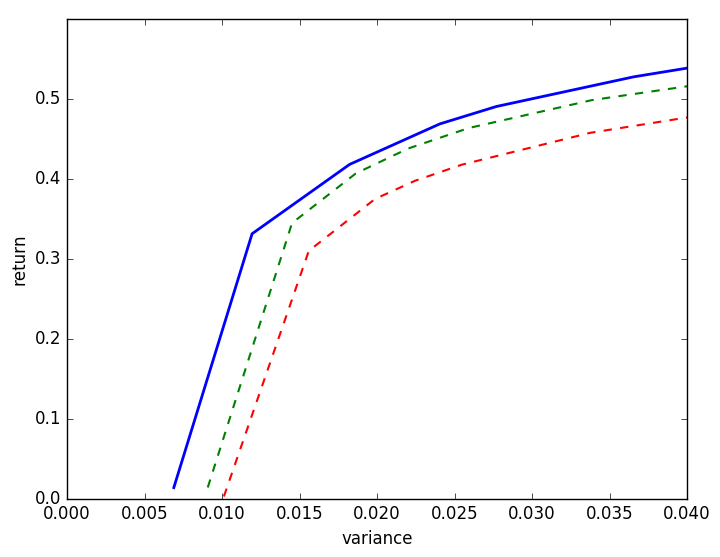}
\caption{Variance v. return for unconstrained (solid line) and cardinality-constrained (dashed lines) Markowitz portfolios as $\gamma$ varies. 
The green dashed uses a model that excludes industrial goods and allows at most 2 stocks from other sectors. 
The red dashed line uses a model that excludes industrial goods and services, and allows at most 2 stocks from other sectors. 
}
\label{fig:mvar}
\end{figure}

Figure \ref{fig:mvar} shows the plot of portfolio variance vs. portfolio return in a Markowitz model as $\gamma$ varies from 0 to 1.5. The blue line shows the Pareto curve 
without cardinality constraints; and the dashed lines show the effects of the constrains on the frontier. 
In particular, the green dashed line is the Pareto frontier using a model that excludes industrial goods and allows at most 2 stocks from other sectors. 
The red dashed line is the Pareto frontier using a model that excludes industrial goods and services, and allows at most 2 stocks from other sectors. 
Progressively more stringent constraints shift us to less efficient regimes; however, the loss in efficiency relative to an unconstrained model is minor (particularly for the first model) and 
can be quantified using the Pareto approach. 

\begin{figure}
\centering
\includegraphics[scale = .52]{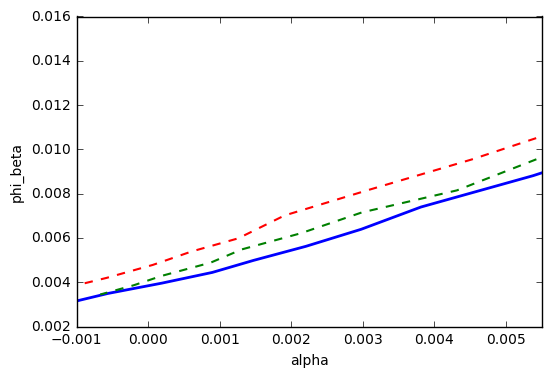}
\caption{$\beta$-VaR ($\alpha$) v. $\beta$-CVaR ($\phi_\beta$) for unconstrained (solid line) and cardinality constrained (dashed lines) CVaR portfolios as $\beta$ varies. 
The green dashed uses a model that excludes industrial goods and allows at most 2 stocks from other sectors. 
The red dashed line uses a model that excludes industrial goods and services, and allows at most 2 stocks from other sectors. }
\label{fig:cvar}
\end{figure}

Figure \ref{fig:cvar} shows the plot of $\beta$-VaR and $\beta$-CVaR values as $\beta$ varies from 0.5 to 0.95. $\beta$-VaR corresponds to $\alpha$ value in the objective and $\beta$-CVaR corresponds to $\phi_\beta = F_\beta(\cdot, \alpha_\beta)$ \cite{rockafellar2002conditional} where
\[ F_\beta(w,\alpha) =  \alpha + \frac{1}{N(1-\beta)}\sum_{j=1}^N [-w^Tr_j - \alpha]_+ .\] 
{Recall that $\alpha$ is the value such that $P(\text{loss} > \alpha) \leq 1-\beta$ and $\phi_\beta$ is the expected loss given loss $\geq \alpha$.}
In Figure \ref{fig:cvar} blue line plots the relation without any constraints; green and red dashed lines correspond to the same constraints as for Figure~\ref{fig:mvar}.
When constraints are imposed, at a given $\alpha$ the value of $\phi_\beta$ increases, meaning that the expected tail loss increases. 
Analogously to the Markowitz case, the increase in tail-loss is small relative what is possible in the unconstrained case, 
and can be quantified using this approach. These results allow fund manages and investors to incorporate realistic 
constraints transaction costs, minimum lot sizes, execution efficiency, and investor preferences directly into portfolio optimization, 
and to evaluate the resulting portfolios against idealized settings.  

\section{Discussion}
In this paper, we proposed a new approach for cardinality-constrained optimization, 
using a relaxation formulation together with modern algorithms for nonconvex optimization. 
The approach is computationally efficient, and guarantees to converge to a stationary point of the original problem. Numerical experiments show that, for small simple cases, the approach can effectively locate global minima for cardinality-constrained Markowitz portfolios. Moreover, in more general cases, the quality of the solutions is high. For larger models, the performance of the approach can be quantified using empirical Pareto frontiers, which show merely a fairly small loss in efficiency compared to the simpler unconstrained portfolios in both Markowitz and CVaR settings. All these suggest high practical applicability of our model and algorithm.  



\section*{APPENDIX}

\section*{Preliminaries from Convex and Variational Analysis}

\begin{definition}[Frechet subdifferential] 
The Frechet subdifferential of a nonconvex function $g$ at $\bar{w}$ is
\[ \hat{\partial} g(\bar{w}) = \{ v: \lim_{u \neq \bar{w}}\inf_{u \to \bar{w}} \frac{g(u) - g(\bar{w}) - \langle v, u-\bar{w}\rangle}{\|u-\bar{w}\|} \geq 0\}.\]
\end{definition}

\begin{definition}[Limiting subdifferential]
The limiting subdifferential of a nonconvex function $g$ at $\bar{w}$ is defined via a closure process
\[ \partial g(\bar{w}) = \{ v: \exists w^j\to \bar{w}, g(w^j) \to g(\bar{w})\]
and  
\[ v^j \in \hat{\partial} g(w^j) \to v \text{ as } k \to \infty.\]
\end{definition}

\begin{definition}[Monotone operator]
A (multivalued) operator $T: \mathbb{R}^n \rightrightarrows\mathbb{R}^n$ is {\it monotone} 
if $\langle u-v, x - y\rangle \geq 0$ for all $u \in Tx$ and $v \in Ty$.
\end{definition}

\begin{proposition}[Monotonicity of subdifferential]
If $g: \R^n \to \bar{\R}$ is closed and convex, then $\partial g$ is monotone.
\end{proposition}

\begin{definition}[Resolvent]
For any operator $T$ and $\alpha > 0$, the resolvent is $(I + \alpha T)^{-1}$. In particular, prox$_{\alpha g} = (I + \alpha \partial g)^{-1}$.
\end{definition}

\begin{theorem}
\label{thm:mono}
If $T$ is monotone, then $(I + \alpha T)^{-1}$ is monotone, single valued and 1-Lipschitz continuous on its domain.
\end{theorem}

\section*{Proof of Lemma \ref{lm:proj}}

The problem can be stated as
\begin{align*}
 &\min_{z_K \in \mathcal{C},K \in \mathcal{K}} \frac{1}{2}\sum_{j \in K}(y_j-z_j)^2 + \frac{1}{2}\sum_{j \in I-K}y_j^2 \\
 \Leftrightarrow &\min_{z_K \in \mathcal{C},K \in \mathcal{K}} \frac{1}{2}\|y_K-z_K\|^2 + \frac{1}{2}\|y_{I-K}\|^2 \\
 \Leftrightarrow &\min_{z_K \in \mathcal{C},K \in \mathcal{K}} \frac{1}{2}\|y_K-z_K\|^2 - \frac{1}{2}\|y_K\|^2 + \frac{1}{2}\|y\|^2.
 \end{align*}

Note that the last term $\frac{1}{2}\|y\|^2$ does not depend on $z_K$, so we can focus on the first two terms, i.e.
$$\min_{z_K \in \mathcal{C},K \in \mathcal{K}}\frac{1}{2}\|y_K-z_K\|^2- \frac{1}{2}\|y_K\|^2.$$

Suppose there is some $K'$ that is different from $K_{opt}$ and denote the corresponding $y$ as $y_{K'}$. Define $f(y)$ and $g(t)$ by 
\[ 
\begin{aligned}
f(y) & = -\frac{1}{2}\|y\|^2 + \min_{z \in \mathcal{C}}\frac{1}{2}\|y-z\|^2, \\
g(t) & = f((1-t)y_{K_{opt}} + t y_{K'}).
\end{aligned}
\]
Then we have 
\[
\begin{aligned}
&f(y_{K'}) - f(y_{K_{opt}}) = g(1) - g(0) = \int_0^1 g'(t) dt, \\
 &g'(t)  = \nabla f((1-t)y_{K_{opt}} + t y_{K'})^T(-y_{K_{opt}} + y_{K'}), 
 \end{aligned}\]
 where
 $\nabla f(y)  = -y + y - z^* = -z^* \in -\mathcal{C}$ given that $\mathcal{C}$ is convex. Since $\mathcal{C} \subset \R_+$, $\nabla f(y)$ is nonpositive in all components. Therefore, $ \nabla f((1-t)v_{K_{opt}} + t v_{K'}) \leq 0$. Further $-v_{K_{opt}} + v_{K'} \leq 0$ because $v_{K_{opt}}$ contains the $k$-largest components of $v$. As a result,
\[ g'(t) \geq 0 \Rightarrow \int_0^1 g'(t) dt \geq 0 \Rightarrow f(v_{K'}) \geq f(v_{K_{opt}}).\]
This shows that $K_{opt}$ must be the optimal choice. 

\addtolength{\textheight}{-10cm}   

\section*{Proof of Lemma~\ref{lem:box} }
We prove each case. 

If $1^Tu_{1:m} < p$, then
\[ 1^Tz^* \geq p > 1^Tu_{1:m} \geq 1^Tu \Rightarrow 1^T(z^*-u) > 0,\]
which implies that the vector $z^* - u$ must have some positive entries. Let $J$ be the index set of those entries. If $1^Tz^* = p + \epsilon > p$, we can decrease $z^*_j-u_j$ for some $j \in J$ to $\max(0, z_j^*-\epsilon-u_j)$ by reducing the value of $z_j^*$. Thus the objective value is decreased, indicating that such a $z^*$ is not optimal. In other words $1^Tz^* = p$ must hold. Hence the problem reduces to 
\[ \min_{z \geq 0, 1^T z=p }\|u-z\|^2 \Leftrightarrow \min_{z \in \Delta_{p}}\|u-z\|^2.\]

If $p \leq 1^Tu_{1:m} \leq q$,
\[ \min_{z \geq 0, p\leq 1^Tz \leq q} \|u-z\|^2 \geq \min_{z \geq 0} \|u-z\|^2 = \|u_{m+1:k}\|^2\]
where the equality can be achieved when $z^*_{1:m} = u_{1:m}$ and $z^*_{m+1:k} = 0$.

If $ 1^Tu_{1:m} > q$, then 
\[ 1^Tu_{1:m} > 1^Tz^* \geq 1^Tz^*_{1:m} \Rightarrow u_{1:m}-z^*_{1:m} > 0,\]
which means that $u_{1:m}-z^*_{1:m}$ must have some positive entries. Let $J$ be the index set of those entries. If $1^Tz^* = q-\epsilon < q$, we can always decrease $u_j - z^*_j$ for some $j \in J$ to $\max( 0, u_j-(z^*_j +\epsilon))$ by increasing the value of $z_j^*$. This indicates $z^*$ is not optimal, hence $1^Tz^* = q$ must hold. Then the problem reduces to
\[ \min_{z \geq 0, 1^Tz = q}\|u-z\|^2 \Leftrightarrow \min_{z \in \Delta_q}\|u-z\|^2.\]

\section*{Proof of Lemma~\ref{lm:station}}
We first need to determine $\partial \delta_C(w)$, by looking at the Frechet subdifferential. The interesting case is when $u \in C$, in which case the expression reduces to
\[  \lim_{u \neq \bar{w}}\inf_{u \to \bar{w}} \frac{ -\langle v, u-\bar{w}\rangle}{\|u-\bar{w}\|} \geq 0, u \in C.\]
Suppose $v \in \hat{\partial} \delta_C(\bar{w})$ and $u \in B$, then 
\begin{align*}
 v \in \hat{\partial} \delta_C(\bar{w})& \Rightarrow \lim \inf_{u \in B \to \bar{w}} -\langle v ,u - \bar{w}\rangle \geq 0 \\
& \Rightarrow  -\langle v ,u - \bar{w}\rangle \geq 0~ \forall u \in B \Rightarrow v \in N_B(\bar{w}).
\end{align*}
Hence $\hat{\partial} \delta_C(w) \subseteq \partial \delta_B(w)$ always holds. 

If further $w \in \text{int} \Delta_1^k$, then as $u \to w$, $u \in B$. Hence 
\begin{align*}
v \in N_B(w)& \Rightarrow \langle v , u -w \rangle \leq 0 ~\forall u \in B \\
 &\Rightarrow \lim_{u \neq w}\inf_{u \to w} \frac{ - \langle v, u-w\rangle}{\|u-w\|} \geq 0 \\
 & \Rightarrow v \in \hat{\partial} \delta_C(w).
 \end{align*}
In other words, $ \partial \delta_B(w) = \hat{\partial} \delta_C(\bar{w}) = \partial \delta_C(\bar{w})$.

\section*{Proof of Lemma~\ref{lem:proxcvx}}

The preliminaries (monotonicity of the subdifferential, definition of resolvent, and Theorem~\ref{thm:mono})  imply that the proximal operator for convex function is single valued. 
We can argue that stationary points and fixed points are equivalent in the  convex case. The arguments below are standard, see e.g. \cite[Proposition 3.1]{combettes2005signal}.
\begin{align*}
 0 \in \nabla g(x) + \partial f(x) \Leftrightarrow & -\nabla f(x) \in \partial g(x) \\
 \Leftrightarrow & (x - \alpha \nabla f(x))  \in x + \alpha \partial g(x) \\
  \Leftrightarrow & (x - \alpha \nabla f(x))  \in (I + \alpha \partial g) (x) \\
 \Leftrightarrow & x = \text{prox}_{\alpha g}(x - \alpha \nabla f(x))
 \end{align*}
 since $(I + \alpha \partial g)^{-1}$ is single valued. 
This proof would not work more generally (in the non-convex case), since the proximal operator may not be single valued.

%

\bibliographystyle{unsrt}
\bibliography{ref}


\end{document}